\documentclass[article,12pt]{amsart}

\usepackage{amsfonts}
\usepackage{xcolor}
\usepackage{amsmath}
\usepackage[margin=2cm]{geometry}

\newtheorem{thm}{Theorem}[section]

\newcommand{\dE}{\mathbb{E}}

\newcommand{\ind}{\mbox{1}\kern-.25em \mbox{I}}
\font\calcal=cmsy10 scaled\magstep1
\def\build#1_#2^#3{\mathrel{\mathop{\kern 0pt#1}\limits_{#2}^{#3}}}
\def\liml{\build{\longrightarrow}_{}^{{\mbox{\calcal L}}}}

\def\videbox{\mathbin{\vbox{\hrule\hbox{\vrule height1ex \kern.5em
\vrule height1ex}\hrule}}}

\email{cristian.coletti@ufabc.edu.br}
\email{manuel.gonzaleznavarrete@ufrontera.cl}
\email{vvazquez@fcfm.buap.mx}

\usepackage{cancel}
\usepackage{verbatim}   

\begin{document}

\title[Asymptotics of a MRW]{Complementary asymptotic analysis for a minimal random walk.}

\author{Cristian F. Coletti}
\address{Universidade Federal do ABC, Center of Mathematics, Computation and Cognition, Sao Paulo, Brazil.}

\author{Manuel Gonz\'alez-Navarrete}
\address{Universidad de La Frontera, Departamento de Matem\'atica y Estad\'istica, Temuco, Chile.}

\author{V\'ictor Hugo V\'azquez Guevara}
\address{Benem\'erita  Universidad Aut\'onoma de Puebla, Facultad de Ciencias Fisico, Puebla, M\'exico.}








	\maketitle

	
\vspace{10pt}

\begin{abstract}
We discuss a complementary asymptotic analysis of the so called minimal random walk. More precisely, we present a version of the almost sure central limit theorem as well as a generalization of the recently proposed quadratic strong laws. In addition, alternative demonstrations of the functional limit theorems will be supplied based on a P\'olya urn scheme instead of a martingale approach.

\end{abstract}

%
\vspace{2pc}
\noindent{\it Keywords}: Minimal Random Walk, Almost Sure Central Limit Theorem, P\'olya urns.
%
%

\section{Introduction}

In this paper we complement the asymptotic analysis developed in \cite{BVMRW} and \cite{coletti3} for the minimal random walk (MRW) proposed by Harbola, Kumar and Lindenberg \cite{kumar} in $2014$. More precisely, we provide the almost sure central limit theorem (ASCLT) for the diffusive and critical regimes as well as the convergence to even moments of a Gaussian random variable,  that generalizes the quadratic strong laws previously obtained in \cite{BVMRW}. Additionally, we achieve an alternative demonstration for the functional central limit theorems (FCLT) through its relation with an urn of P\'olya type in the spirit of \cite{baur} which was originally proved in \cite{BVMRW}. It must be mentioned that most of our technical display will be based; to some extent, on some side properties provided by \cite{BVMRW}, situation that will be properly treated in Appendix $A$.

The rest of the paper is organized as follows: Section \ref{sectionMRW}  presents to the MRW as well as its basic conditional properties. Section \ref{sectionMR} provides versions of the ASCLT and convergence to even moments of Gaussian distribution in the diffusive and critical regimes. Section \ref{sectionPU} explores the relation of the MRW and a P\'olya urn scheme with balls of two colors to provide alternative proofs of the versions of the FCLT obtained in \cite{BVMRW}. Four appendices are considered: Appendix A summarizes some properties taken from \cite{BVMRW} and   provides the conditional moments of the MRW. Appendices B and C deal with the demonstrations corresponding to results of Section \ref{sectionMR}. Last appendix is devoted to the proof of the FCLT's of Section \ref{sectionPU}.



\section{The Minimal Random Walk} \label{sectionMRW}
In this section we present the main topic of this paper: The MRW, whose movement and position at time $n\geq 1$ are denoted by $X_n$ and $S_n$ respectively behaves as follow: 

\begin{itemize}
\item At time $n=0$ a particle following the MRW is at the origin; i.e., $S_0=0$.
\item At time $n=1$, the particle decides to move one unit right with probability $s\in(0,1)$ or remain at the origin with probability $1-s$. Then $S_1=X_1$ has the Bernoulli distribution with parameter $s$. 
\item At any time $n\geq 2$, the particle chooses uniformly at random some point in the past; let us say $1\leq k(n) \leq n-1$. Then, its following movement; $X_{n}$, will be conditioned by the observed value of $X_{k(n)}$ as follows:

On the one hand, if $X_{k(n)}=0$ then $X_n$ posses a Bernoulli conditional distribution with parameter $q\in(0,1]$ (It should be underlined that the case $q = 0$ will not be treated in this paper). While, on the other hand, if $X_{k(n)}=1$ then, the conditional distribution of $X_n$ is also Bernoulli but with parameter $p\in[0,1]$, where $p$ does not necessarily is equal to $q$. 

\item The position of the particle at time $n\geq 1$ satisfies that
\begin{equation}\label{position}
S_{n}=S_{n-1}+X_{n},
\end{equation}%
with $S_0=0$ and $X_n\in\left\{0,1\right\}$.
\end{itemize}






In order to consider the martingale approach of \cite{BVMRW} for the MRW, let $\left(\mathcal{F}_n\right)$ be the increasing sequence of $\sigma$-algebras $\mathcal{F}_n=\sigma\left(X_1,\ldots,X_n\right)$; i.e. $\mathcal{F}_n$ represents the knowledge up to time $n$.

From \cite{kumar}, we know that for $x=0,1$ , $\mathbb{P}\left(X_{n+1}=x|\mathcal{F}_n\right)=1-x+(2x-1)\left(q+\alpha\frac{S_n}{n}\right)$, then it may be found that the conditional expected movement and position of the particle at any time are; respectively, given by

\begin{equation}\label{expmov}
\mathbb{E}\left[X_{n+1}|\mathcal{F}_n\right]=q+\alpha\frac{S_n}{n} \hspace{.5cm}\text{a.s.}
\end{equation}
and
\begin{equation}\label{expos}
\mathbb{E}\left[S_{n+1}|\mathcal{F}_n\right]=q+\gamma_n S_n \hspace{.5cm}\text{a.s.}
\end{equation}%




where, $\alpha:=p-q$ and $\gamma_n=\frac{n+\alpha}{n}$.
\\

Let us now consider the discrete time scalar martingale $(M_n)$, given for $n\geq 0$, by 

\begin{equation} 
\label{martingale}
M_n=a_n S_n-q A_n,
\end{equation}
where $a_1=1$ and, for $n\geq 2$

\begin{equation}\label{an}
a_n=\prod_{k=1}^{n-1}\gamma_k^{-1}=\frac{\Gamma(n)\Gamma(\alpha+1)}{\Gamma(n+\alpha)},
\end{equation}%
where $\Gamma$ stands for the Euler gamma function. In addition, $A_0=0$ and, for $n\geq 1$,
\begin{equation}
\label{bigAn}
A_n=\sum_{k=1}^n a_k.
\end{equation}

The study of asymptotic properties of the scalar MRW is usually done in 3 regimes:
\begin{enumerate}
\item If $\alpha<1/2$ then we are in the diffusive regime.
\item If $\alpha=1/2$ then the regime is the critical one, and
\item If $\alpha>1/2$ then the regime is called superdiffusive.
\end{enumerate}


Provided results in following sections will be; basically of the almost sure central limit type, convergence to even moments of Gaussian distribution and functional central limit theorems.

\section{Deeper results via the martingale approach} \label{sectionMR}

In the paper \cite{coletti3}, it has been proved versions of the central limit theorem and the law of iterated logarithm in the case $\alpha \le 1/2$, which were complemented by \cite{BVMRW}, including law of large numbers, functional central limit theorems, an asymptotic analysis of the center of mass and a better understanding of the superdifussive case $\alpha > 1/2$ by providing the almost sure and in $L^2$ convergence as well as the Gaussian fluctuations around its limit, in addition the first two moments of this random variable were explicitly calculated. Furthermore, in \cite{miyazaki} it has been considered the case $q=0$ where a central limit theorem and a law of iterated logarithm were obtained.

In this section, we will display proper versions of the almost sure central limit theorem and of a generalization of the recently proposed quadratic strong laws under the condition $\alpha \leq 1/2$. For this, in all the sequel we will make use of the asymptotic variance
\begin{equation}
\label{varianza}
\sigma^2:= \frac{q(1-p)}{(1-\alpha)^2}
\end{equation}

\subsection{The diffusive regime}

We initiate the asymptotic analysis by dealing with the almost sure central limit theorem (ASCLT) for the diffusive regime $\left(\alpha <1/2\right)$.
\\

\begin{thm}\label{CVG1}
If $\alpha<1/2$ then, we have the following almost sure convergence of measures
\begin{equation}\label{res1}
\frac{1}{\log n}\sum_{k=1}^n \frac{1}{k}\delta_{\sqrt{k}\left(\frac{S_k}{k}-\frac{q}{1-\alpha}\right)}\Longrightarrow N\left(0,\frac{\sigma^2}{1-2\alpha}\right)
\end{equation}%
where $\sigma^2$ is given in \eqref{varianza}.
\end{thm}

In a complementary way, we provide the corresponding convergence of moments, which generalizes the quadratic strong law found in \cite{BVMRW}.

\begin{thm}\label{CVG2}
If $\alpha <1/2$, then we have; for $r\geq 1$, the following almost sure convergence 

\begin{equation}\label{res2}
\lim_{n\rightarrow \infty} \frac{1}{\log n}\sum_{k=1}^n k^{r-1} \left(\frac{S_k}{k}-\frac{q}{1-\alpha}\right)^{2r}=\frac{\sigma^{2r}(2r)!}{2^r r! (1-2\alpha)^{r} }
\end{equation}%
\end{thm}

\subsection{The critical regime}

We will investigate now, analogue results to those presented in previous subsection for the critical regime, starting with the ASCLT.

\begin{thm}\label{CVG11}
If $\alpha = 1/2$, then we have that almost surely 
\begin{equation}\label{as11}
\frac{1}{\log \log n}\sum_{k=2}^n \frac{1}{k\log k}\delta_{\sqrt{\frac{k}{\log k}} \left( \frac{S_k}{k}-2q\right)}\Longrightarrow N\left(0,4q(1-p)\right)\end{equation}
\end{thm}

Next theorem establishes the convergence to  even moments of Gaussian distribution

\begin{thm}\label{CVG21}
If $\alpha = 1/2$, then the following almost sure convergence holds
\begin{equation}\label{qsl11}
\lim_{n\rightarrow \infty} \frac{1}{\log \log n}\sum_{k=2}^n \left(\frac{1}{\log k}\right)^{r+1} k^{r-1} \left( \frac{S_k}{k}-2q\right)^{2r}=\frac{(4q(1-p))^{r}(2r)!}{2^r r!}
\end{equation}
where $r\geq 1$.
\end{thm}

\section{Relation of the MRW with P\'olya urns} \label{sectionPU}

In this section, we will perform a counterpart display of the ideas exposed in \cite{baur} for the elephant random walk. Hence, in order to establish a relationship between the MRW and an urn of P\'olya type, let us consider 

\begin{itemize}
\item An urn with balls of two colors:  red and blue.
\item At $n=1$ we put in a blue one with probability (w.p.) $s$.
\item At $n>1$ we choose uniformly at random one ball, observe its color, put it back and if it is red, then we put another red w.p. $1-q$.
\item If it is blue then we put another blue w.p. $p$.
\end{itemize}

Hence, $S_n$ is the number of blue balls after $n$ extractions which evolves as the MRW. Furthermore, we obtain the so-called mean replacement matrix for the earlier urn associated with the MRW:

\begin{equation}
A:=
\begin{pmatrix}
1-q & 1-p\\
q & p
\end{pmatrix}
\end{equation}
whose eigenvalues are $\lambda_1:=1$ and $\lambda_2:=\alpha$ with $\lambda_2 \leq \lambda_1$.

In succeeding subsection, the theory for generalized P\'olya urns established in \cite{janson} will be employed in order to find alternative demonstrations for the functional central limit theorems in both the diffusive and critical regimes as well as the almost sure convergence in the superdiffusive regime. In all the sequel $D([0,\infty[)$ denotes to the Skorokhod space of right-continuous functions with left-hand limits.

\subsection{Limit theorems via the P\'olya urns scheme}

Subsequent results were obtained in \cite{BVMRW} via the martingale approach considered in sections \ref{sectionMRW} and \ref{sectionMR}, nevertheless in Appendix $D$ the link between the MRW and the previously urn scheme will be employed in order to prove them.

\begin{thm} \label{fclt1}
If $\alpha<1/2$ then, we have the distributional convergence in $D([0,\infty[)$,

$$
\left( \sqrt{n}\Big(\frac{S_{\lfloor nt \rfloor}}{\lfloor nt \rfloor}-\frac{q}{1-a}\Big), t \geq 0\right) \Longrightarrow \big( W_t, t \geq 0 \big)
$$

where $\big( W_t, t \geq 0 \big)$ is a real-valued centered Gaussian process starting at the origin with covariance given, for all $0<s \leq t$, by

$$
\dE[W_s W_t]= \frac{\sigma^2}{(1-2a)t} \Big(\frac{t}{s}\Bigr)^a.
$$

In particular, we have the asymptotic normality

$$ \sqrt{n}\Big(\frac{S_n}{n}-\frac{q}{1-a}\Big) \liml N\Big(0,\frac{\sigma^2}{1-2a}\Big)$$

\end{thm}

\begin{thm} \label{fclt2}

If $\alpha=1/2$ then, we have the distributional convergence in $D([0,\infty[)$,

$$ \left( \sqrt{\frac{n^t}{\log n}}\Big(\frac{S_{\lfloor n^t \rfloor}}{\lfloor n^t \rfloor}-2q\Big), t \geq 0\right) \Longrightarrow \big( 2\sqrt{q(1-p)} B_t, t \geq 0 \big)$$

where $\big( B_t, t \geq 0 \big)$ is a standard Brownian motion.
\\
In particular, we have the central limit theorem

$$ \sqrt{\frac{n}{\log n}}\Big(\frac{S_n}{n}-2q\Big) \liml N\big(0,4q(1-p)\big)$$

\end{thm}


\begin{thm} \label{fclt3}

If $\alpha>1/2$ then, we have the almost sure convergence,

$$ n^{1-a}\Big(\frac{S_{n}}{n}-\frac{q}{1-a}\Big) \Longrightarrow L
$$

where  $L$ is a non-degenerated random variable and whose two first moments are calculated in Theorem $2.6$ of \cite{BVMRW}.

\end{thm}

\section*{Acknowledgments}
Cristian F. Coletti was supported by grant \#2017/10555-0 S$\tilde{a}$o Paulo Research Foundation (FAPESP). Manuel Gonz\'alez-Navarrete was partially supported by Fondecyt Iniciaci\'on 11200500.

\section*{Appendix A \\ Initial Analysis}
\renewcommand{\thesection}{\Alph{section}}
\renewcommand{\theequation}{\thesection.\arabic{equation}}
\setcounter{section}{1}
\setcounter{equation}{0}

In this Appendix, we will state some useful facts on the MRW in order to achieve Theorems \ref{CVG1} - \ref{CVG21}: 

For each $n\geq 1$ and $k\geq 1$, the definition of the MRW ensures that, $X_n \leq 1$, then $0\leq S_n \leq n$ and $X_n=X_n^k$. Moreover,
$$\mathbb{E}\left[ X_{n+1}^k |\mathcal{F}_n\right]=q+\alpha\frac{S_n}{n} \hspace{.5cm}\text{a.s.} $$
and
$$\mathbb{E} \left[S_{n+1}^{2k}|\mathcal{F}_n\right]=S_n^{2k}+\left[ q+\alpha\frac{S_n}{n} \right] \sum_{j=1}^{2k} \binom{2k}{j} S_n^{2k-j}. $$ 

In addition, the martingale $(M_n)$ may be rewritten as $M_n=\sum_{k=1}^n a_k \varepsilon_k$, where $\varepsilon_{n+1}=X_{n+1}-\left(q+\alpha\frac{S_n}{n} \right).$ It also holds that 
$$\mathbb{E}\left[\varepsilon_{n+1}^{k}|\mathcal{F}_n\right]=\sum_{j=0}^{k-2} \binom{k}{j} \left(q+\alpha\frac{S_n}{n} \right)^{j+1} (-1)^j+(-1)^{k-1}\left(q+\alpha\frac{S_n}{n} \right)^k (k-1)<\infty.$$
Therefore, Theorems $2.1$ and $2.3$ of \cite{BVMRW} guide us to
\begin{equation} \label{segundo}
\lim_{n\rightarrow \infty} \mathbb{E}\left[\varepsilon^2_{n+1}|\mathcal{F}_n\right] =\sigma^2
\hspace{1cm} \text{a.s.},
\end{equation}
finally, a very useful consequence of Lemma $B.1$ of \cite{ERWBercu} which will be repeatedly used is 
\begin{equation} \label{relaan}
\left| \frac{A_n}{n a_n}-\frac{1}{1-\alpha} \right|\sim \frac{1}{n^{1-\alpha}}.
\end{equation}

In addition, by continuing the notation of \cite{BVMRW}, let us consider the following entities for each $n \geq 1$:
$$v_n:=\sum_{k=1}^n a_k^2.$$
and the explosion coefficient associated with martingale $(M_n)$
$$f_n:=a_n^2/v_n$$.

\section*{Appendix B \\ASCLT for the diffusive regime}
\renewcommand{\thesection}{\Alph{section}}
\renewcommand{\theequation}{\thesection.\arabic{equation}}
\setcounter{section}{2}
\setcounter{equation}{0}


{ \bfseries Proof of theorem \ref{CVG1} and \ref{CVG2}}:
For this regime, it has been previously obtained \cite{BVMRW} that, as $n\rightarrow \infty$:
\begin{equation}\label{vndif}
\frac{v_n}{n^{1-2\alpha}}\rightarrow \ell:=\frac{\Gamma^2(\alpha+1)}{1-2\alpha}
\end{equation}
and
\begin{equation}\label{relfn}
f_n \sim \frac{1-2\alpha}{n}\rightarrow 0.
\end{equation}
Additionally, \eqref{relaan} implies the following useful relation:
\begin{equation}\label{ap11}
\frac{M_k}{\sqrt{v_{k-1}}} \sim \sqrt{\frac{1-2\alpha}{k}} \left(S_k - k\frac{q}{1-\alpha}\right)
\end{equation}

In order to utilize Theorem $A.1$ of \cite{yo} for proving Theorem \ref{CVG1}, we note from last statements that for $\eta >0$

     \begin{eqnarray}
    \displaystyle\displaystyle\sum_{k=1}^{\infty}\frac{1}{v_k}\dE\left[\vert\Delta M_{k}\vert^2 \mathbb{I}_{\vert\Delta M_{k}\vert\geq\eta\sqrt{v_k}}\vert\mathcal{F}_{k-1}\right]       \leq \frac{1}{\eta^2}\displaystyle\sum_{k=1}^{\infty}\frac{1}{v_k^2}\dE\left[\vert\Delta M_{k}\vert^4 \vert\mathcal{F}_{k-1}\right] \\[0.5cm]
           \leq \displaystyle\sup_{k\geq 1}\dE\left[\varepsilon^4_{k} \vert\mathcal{F}_{k-1}\right]\frac{1}{\eta^2}\displaystyle\sum_{k=1}^{\infty} \frac{a_k^4}{v_k^2 }  \sim \frac{\displaystyle\sup_{k\geq 1}\dE\left[\varepsilon^4_{k} \vert\mathcal{F}_{k-1}\right]}{\eta^2}\displaystyle\sum_{k=1}^{\infty} \frac{(1-2\alpha)^2}{k^2} < \infty.
     \end{eqnarray}
On the same fashion, it may be found that
$$
 \displaystyle\displaystyle\sum_{k=1}^{\infty}\frac{1}{v_k^2}\dE\left[\vert\Delta M_{k}\vert^4 \mathbb{I}_{\vert\Delta M_{k}\vert\leq\sqrt{v_k}}\vert\mathcal{F}_{k-1}\right]       \leq \displaystyle\sum_{k=1}^{\infty}\frac{a_k^4}{v_k^2}\dE\left[\varepsilon^4_{k} \vert\mathcal{F}_{k-1}\right]<\infty.$$

Hence, Theorem $A.1$ of \cite{yo} let us conclude that

\begin{equation}
  \displaystyle\frac{1}{\log v_{n}}\sum_{k=1}^n f_k \delta_{M_k/\sqrt{v_{k-1}}} \Rightarrow G  \ \  \text{ a.s}
\end{equation}
where $G$ is the Gaussian distribution with mean zero and variance $\sigma^2$. Finally, convergence \eqref{vndif}, \eqref{relfn} and relation \eqref{ap11}  imply 

\begin{equation}
  \displaystyle\frac{1}{\cancel{(1-2\alpha)}\log n}\sum_{k=1}^n \frac{\cancel{(1-2\alpha)}}{k} \delta_{\sqrt{1-2\alpha}\left(\frac{S_k}{k}-\frac{q}{1-\alpha}\right)\sqrt{k}} \Rightarrow G  \ \  \text{ a.s}
\end{equation}
which conduces us to convergence \eqref{res1}.
\\

For obtaining Theorem \ref{CVG2}, we may observe from \eqref{segundo} and \eqref{relaan}, together with Theorem $3$ of \cite{cvgm} (or  Theorem $A.2$ of \cite{yo}) that; for $r\geq 1$:

\begin{equation}
\frac{1}{\log n}\sum_{k=1}^n f_k \left(\frac{M_k^2}{v_{k-1}}\right)^r \rightarrow \frac{\sigma^{2r}(2r)!}{2^r r!,}
\end{equation}
however, definition of $M_n$ and $f_n$, convergences \eqref{vndif} and \eqref{relfn} lead us to

\begin{equation} 
\frac{1}{\log n}\sum_{k=1}^n \frac{1}{k^{r+1}} \left( S_k-q\frac{A_k}{a_k} \right) ^{2r} \rightarrow \frac{\sigma^{2r}(2r)!(1-2\alpha)}{2^r r! (1-2\alpha)^{r+1}},
\end{equation}
that implies

\begin{equation} \label{cvgmdif1}
\frac{1}{\log n}\sum_{k=1}^n k^{r-1} \left( \frac{S_k}{k}-q\frac{A_k}{k a_k} \right) ^{2r} \rightarrow \frac{\sigma^{2r}(2r)!}{2^r r! (1-2\alpha)^{r}}.
\end{equation}

Nevertheless, we have that
\begin{equation} \label{cvgmdif2}
\frac{1}{\log n}\sum_{k=1}^n k^{r-1} \left( \frac{S_k}{k}-\frac{q}{1-\alpha} \right) ^{2r} = \frac{1}{\log n}\sum_{k=1}^n k^{r-1} \sum_{i=0} ^{2r} \binom{2r}{i} \left( \frac{S_k}{k}-q\frac{A_k}{k a_k} \right) ^{i} \left( q\frac{A_k}{k a_k}-\frac{q}{1-\alpha} \right) ^{2r-i}.
\end{equation}

However, it may be noticed from \eqref{relaan}; for $0\leq i < 2r$, that
\begin{eqnarray*}
\frac{1}{\log n} \sum_{k=1} ^{n} k^{r-1} \binom{2r}{i} \left( \frac{S_k}{k}-q\frac{A_k}{k a_k} \right) ^{i} \left( q\frac{A_k}{k a_k}-\frac{q}{1-\alpha} \right) ^{2r-i} \\
\leq \frac{C}{\log n} \sqrt{\sum_{k=1} ^{n} k^{i-1} \left( \frac{S_k}{k}-q\frac{A_k}{k a_k} \right) ^{2i}} \sqrt{\sum_{k=1} ^{n} k^{2r-i-1} \left( q\frac{A_k}{k a_k}-\frac{q}{1-\alpha} \right) ^{4r-2i}}\\
\sim \frac{C}{\log n} \sqrt{\sum_{k=1} ^{n} k^{i-1} \left( \frac{S_k}{k}-q\frac{A_k}{k a_k} \right) ^{2i}} \sqrt{\sum_{k=1} ^{n} k^{2r-i-1} \left( \frac{1}{k^{2(1-\alpha)}} \right) ^{2r-i}} \\
\rightarrow 0,
\end{eqnarray*}%
since $2(1-\alpha)(2r-1)-(2r-i-1)>1$. Hence, we infer that

\begin{eqnarray*}
\lim_{n\rightarrow \infty}\frac{1}{\log n}\sum_{k=1}^n k^{r-1} \left( \frac{S_k}{k}-\frac{q}{1-\alpha} \right) ^{2r}\\
=\lim_{n\rightarrow \infty}\frac{1}{\log n}\sum_{k=1}^n k^{r-1} \left( \frac{S_k}{k}-q\frac{A_k}{k a_k} \right) ^{2r} \\
=\frac{\sigma^{2r}(2r)!}{2^r r! (1-2\alpha)^{r}}.
\end{eqnarray*}

\section*{Appendix C \\ ASCLT for the critical regime}
\renewcommand{\thesection}{\Alph{section}}
\renewcommand{\theequation}{\thesection.\arabic{equation}}
\setcounter{section}{3}
\setcounter{equation}{0}
{ \bfseries Proof of theorems \ref{CVG11} and \ref{CVG21}}:
On this regime, it is useful to highlight that \cite{BVMRW}
\begin{equation} \label{vncrit}
\frac{v_n}{\log n}\rightarrow \frac{\pi}{4},
\end{equation}
which leads us to ensure that
\begin{equation}\label{relfncrit}
 f_n \sim \frac{1}{n \log n}\rightarrow 0,
\end{equation}
 and also that $\sum_{k=1}^\infty \frac{a_k^2}{v_n^2}<\infty$. Therefore, from analogue arguments to those exposed in Appendix $B$ we conclude that hypothesis of the ASCLT expressed in Theorem $A.1$ of \cite{yo} are fulfilled. Hence, from such ASCLT, relation \eqref{relaan} and convergence \eqref{vncrit} we infer \eqref{as11}.

In order to demonstrate Theorem \ref{CVG21}, we may note; from similar argumentation than in Appendix $B$ , that \eqref{segundo} and \eqref{relaan}, together with Theorem $3$ of \cite{cvgm}  and convergences \eqref{vncrit} and \eqref{relfncrit}; for $r\geq 1$, that
\begin{equation} \label{cvgmdif11}
\frac{1}{\log \log n}\sum_{k=2}^n k^{r-1}\left(\frac{1}{\log k} \right)^{r+1} \left( \frac{S_k}{k}-q\frac{A_k}{k a_k} \right) ^{2r} \rightarrow \frac{\sigma^{2r}(2r)!}{2^r r!}.
\end{equation}
However, we may observe that

\begin{eqnarray} \label{cvgmdif21}
\frac{1}{\log \log n}\sum_{k=2}^n k^{r-1}\left(\frac{1}{\log k} \right)^{r+1} \left( \frac{S_k}{k}-2q \right) ^{2r} \\
= \frac{1}{\log \log n}\sum_{k=2}^n k^{r-1}\left(\frac{1}{\log k} \right)^{r+1} \sum_{i=0} ^{2r} \binom{2r}{i} \left( \frac{S_k}{k}-q\frac{A_k}{k a_k} \right) ^{i} \left( q\frac{A_k}{k a_k}-2q \right) ^{2r-i}. \notag
\end{eqnarray}

For $0\leq i < 2r$, we observe that Cauchy-Schwarz inequality, relation \eqref{relaan} and convergence \eqref{cvgmdif11} imply that
$$ \frac{1}{\log \log n}\sum_{k=2}^n k^{r-1} \left(\frac{1}{\log k} \right)^{r+1} \binom{2r}{i} \left( \frac{S_k}{k}-q\frac{A_k}{k a_k} \right) ^{i} \left( q\frac{A_k}{k a_k}-2q \right) ^{2r-i} \rightarrow 0 $$
because
$$\lim _{n\rightarrow \infty} \sum_{k=1}^n \left( \frac{1}{\log k} \right)^{2r-i+1} \frac{1}{k}<\infty$$
since $0\leq i < 2r$ and $r\geq 1$.
\\

All this leads us to find out that

\begin{eqnarray*}
\lim_{n\rightarrow \infty}\frac{1}{\log \log n}\sum_{k=2}^n k^{r-1}\left(\frac{1}{\log k} \right)^{r+1} \left( \frac{S_k}{k}-2q \right) ^{2r}\\
=\lim_{n\rightarrow \infty}\frac{1}{\log \log n}\sum_{k=2}^n k^{r-1}\left(\frac{1}{\log k} \right)^{r+1} \left( \frac{S_k}{k}-q\frac{A_k}{k a_k} \right) ^{2r} \\
=\frac{\sigma^{2r}(2r)!}{2^r r!}.
\end{eqnarray*}

\section*{Appendix D \\ The functional central limit theorem based on a P\'olya urn scheme}
\renewcommand{\thesection}{\Alph{section}}
\renewcommand{\theequation}{\thesection.\arabic{equation}}
\setcounter{section}{4}
\setcounter{equation}{0}
In the sequel, we will employ the notation of \cite{janson} in an attempt to facilitate consultation of the theory contained therein and be able to expand the details of this Appendix; if it is necessary. Hence, at this setting, it is not hard to see that conditions $(A1)-(A6)$ (of that paper) hold for the MRW and his eigenvalues. Additionally, we find out that right and left eigenvectors given by
$$v_1=\frac{q}{1-\alpha}\left(\frac{1-p}{q},1\right)^t, v_2=\frac{1}{2}(-1,1)^t$$
and
$$u_1=(1,1), u_2=\frac{2(p-1)}{\alpha-1}\left(\frac{q}{p-1},1\right);$$

respectively, are such that $u_1\cdot v_1=u_2 \cdot v_2 =1$. Finally, we denote by $W_n=(P_n,S_n)$ to the composition of the urn at stage $n\geq 1$ where $P_n$ states for the number of red balls and $S_n$ is the number of blue ones.

{ \bfseries Proof of theorem \ref{fclt1} (the diffusive regime)}:

For this regime we will make use of Theorem 3.31 i) of \cite{janson} that asserts that 

\begin{equation}\label{th331i}
\left( n^{-1/2}\left( W_{\lfloor nt \rfloor}-nt v_1 \right), t\geq 0 \right) \Longrightarrow \left( V(t), t\geq 0 \right).
\end{equation}
where $\left(V(t),t \geq0 \right)$ is a centered Gaussian vector process such that $V(0)=(0,0)$ and with covariance structured given for $0<s\leq t$ by (see Remark 5.7 of \cite{janson})

\begin{equation}\label{cov1}
\mathbb{E}\left[V(s)V^t(t)\right]=s\Sigma_I e^{\log(t/s)A^t},
\end{equation}
where $\Sigma_I$ was defined in equation $(2.15)$ of \cite{janson} and given for the MRW by

\begin{equation}
\Sigma_I=\frac{q(1-p)}{(1-2\alpha)(1-\alpha)^2}
\begin{pmatrix}
	1 & -1 \\
	-1 & 1 
	\end{pmatrix}.
	\end{equation}
Also, we have that
$$\left( \frac{t}{s}\right) ^{A^t}=P^{-1}\left( \frac{t}{s}\right)^{\alpha} I_2 P$$
where $I_2$ is identity matrix or order $2$ and
\begin{equation}
P=I_2-v_1 u_1^t
=\frac{1}{1-\alpha}
\begin{pmatrix}
q & p-1 \\
-q & 1-p
\end{pmatrix}
.
\end{equation}
Hence, \eqref{cov1} implies
\begin{equation}\label{cov1f}
\mathbb{E}\left[V(s)V^t(t)\right]=s\left( \frac{t}{s}\right)^{\alpha}\frac{q(1-p)}{(1-2\alpha)(1-\alpha)^2}
\begin{pmatrix}
	1 & -1 \\
	-1 & 1 
	\end{pmatrix},
\end{equation}

On the other hand, we observe that $S_n=(0,1)W_n^t$, which together with the continuous mapping theorem leads us to the conclusion of Theorem \ref{fclt1}.
\\

{ \bfseries Proof of theorem \ref{fclt2} (the critical regime)}:

On this regime, Theorem 3.31 ii) of \cite{janson} conduces us to 

$$ \left( \sqrt{\frac{1}{n^t \log n}}\Big(W_{\lfloor n^t \rfloor}-n^t v_1\Big), t \geq 0\right) \Longrightarrow \big( V(t), t \geq 0 \big),$$

where $\left(V(t),t \geq0 \right)$ is a centered Gaussian vector process such that $V(0)=(0,0)$ and with covariance structured given for $0<s\leq t$ by (see equation (3.27) of \cite{janson})

\begin{equation}
\mathbb{E}\left[V(s) V^t(t) \right]=
4sq(1-p)
\begin{pmatrix}
1 & -1\\
-1 & 1
\end{pmatrix}
.
\end{equation}

By considering the same applied transformation as in the diffusive regime, the continuous mapping theorem implies the conclusion of Theorem \ref{fclt2}.

{ \bfseries Proof of Theorem \ref{fclt3} (the superdiffusive regime)}:
This theorem follows essentially from Theorem $3.24$ of \cite{janson}. More precisely, since at this regime the lowest eigenvalue is greater than the half of the largest one (that is $\alpha > 1/2$) ,  we have that, almost surely

$$n^{-\alpha}\left( W_{n}-n v_1 \right)\rightarrow \hat{W},$$

where $\hat{W}$ is a random vector with values in $\left\{w\in \mathbb{R}^2 / w=\lambda(-1,1)\text{ for some } \lambda\in\mathbb{R}-\left\{ 0 \right\} \right\}$, therefore (again) the continuous mapping theorem conduces us to conclusion of Theorem \ref{fclt3}.

\end{document}